\def\lastupdate{15 Sep.\ 2006}
\documentclass[draft]{amsart}
\usepackage{amssymb,amsmath,amscd}
\theoremstyle{plain}
\newtheorem{thm}{Theorem}[section]
\newtheorem{lem}[thm]{Lemma}
\newtheorem{cor}[thm]{Corollary}
\newtheorem*{AT}{The Annihilator Theorem}
\newtheorem*{claim}{Claim}
\theoremstyle{definition}
\newtheorem{dfn}[thm]{Definition}
\numberwithin{equation}{thm}

\newcommand{\ass}{\operatorname{Ass}}
\newcommand{\coker}{\operatorname{Coker}}
\newcommand{\depth}{\operatorname{depth}}
\newcommand{\ext}{\operatorname{Ext}}
\newcommand{\height}{\operatorname{ht}}
\renewcommand{\hom}{\operatorname{Hom}}
\mathchardef\qtn = "103A
\newcommand{\rank}{\operatorname{rank}}
\newcommand{\spec}{\operatorname{Spec}}
\newcommand{\supp}{\operatorname{Supp}}
\begin{document}

\title{On Faltings' annihilator theorem}

\author{Kawasaki, Takesi}

\address{Department of Mathematics and Information Sciences,
  Tokyo Metropolitan University,
  Minami-Ohsawa 1-1,
  Hachioji,
  Tokyo 192-0397,
  Japan}

\email{kawasaki@comp.metro-u.ac.jp}
\urladdr{http://www.comp.metro-u.ac.jp/\textasciitilde kawasaki/}


\dedicatory{Dedicated to Professor Shiro Goto
  on the occasion of his sixtieth birthday}

\date{\lastupdate}

\thanks{This work was supported
  by Japan Society for the Promotion of Science
  (the Grant-in-Aid for Scientific Researches (C)(2)~16540032)}


\begin{abstract}
  In the present article,
    the author shows that Faltings' annihilator theorem holds for
    any Noetherian ring $A$
    if $A$ is universally catenary;
    all the formal fibers of all the localizations of~$A$ are Cohen-Macaulay;
    and the Cohen-Macaulay locus of each finitely generated $A$-algebra
    is open.
\end{abstract}

\maketitle

\section{Introduction}
Throughout the present article,
  $A$ always denotes a commutative Noetherian ring.
We say that the annihilator theorem holds for~$A$
  if it satisfies the following proposition~\cite{MR82f:13003}.

\begin{AT}
  Let $M$ be a finitely generated $A$-module,
    $n$ an integer and
    $Y$, $Z$ subsets of $\spec A$
    which are stable under specialization.
  Then the following statements are equivalent\textup:
    \begin{itemize}
      \item[(1)] $\height \mathfrak p / \mathfrak q + \depth M_{\mathfrak q}
        \geq n$
        for any $\mathfrak q \in \spec A \setminus Y$
        and $\mathfrak p \in V(\mathfrak q) \cap Z$\textup;
      \item[(2)] there is an ideal $\mathfrak b$ in~$A$
        such that $V(\mathfrak b) \subset Y$ and
        $\mathfrak b$ annihilates local cohomology modules
        $H_Z^0(M)$, \dots, $H_Z^{n-1}(M)$.
    \end{itemize}
\end{AT}

Faltings~\cite{MR58:22058} proved that
  the annihilator theorem holds for~$A$
  if $A$ has a dualizing complex or
  if $A$ is a homomorphic image of a regular ring
  and that (2) always implies~(1).
Several authors~\cite{MR2002b:13027,%
  MR99h:13020,%
  MR2052396,%
  MR2229481,%
  MR94e:13033}
  tried to improve Faltings' result.
In this article, the author shows the following
\begin{thm} \label{thm:1.1}
  The annihilator theorem holds for $A$
    if
      \begin{itemize}
        \item[(C1)] $A$ is universally catenary\textup;
        \item[(C2)] all the formal fibers of all the localizations of~$A$
          are Cohen-Macaulay\textup; and
        \item[(C3)] the Cohen-Macaulay locus
          of each finitely generated $A$-algebra
          is open.
      \end{itemize}
\end{thm}

These conditions are not only sufficient
  but also necessary
  for the annihilator theorem.
Indeed, Faltings~\cite{MR82f:13003} showed that $A$ satisfies (C1)--(C3)
  whenever the annihilator theorem holds for
  each essentially of finite type $A$-algebra.

These conditions are also related to the uniform Artin-Rees theorem
  and the uniform Brian\c{c}on-Skoda theorem.
We give an affirmative answer to the conjecture
  of Huneke~\cite[Conjecture~2.13]{MR93b:13027}
  in the last section.

\section{Preliminaries}

First we recall the definition of the local cohomology functor.
A subset $Z$ of $\spec A$ is said to be stable under specialization
  if $\mathfrak p \in Z$ implies $V(\mathfrak p) \subset Z$.
Let $M$ be an $A$-module and
  $Z$ a subset of $\spec A$ which is stable under specialization.
Then we put
    $$
      H_Z^0(M) = \{m \in M \mid \supp Am \subset Z\}.
    $$
It is an $A$-submodule of~$M$ and $H_Z^0(-)$ is a left exact functor.

\begin{dfn}[{\cite[p.~223]{MR36:5145}}]
  The local cohomology functor $H_Z^p(-)$ with respect to~$Z$
    is the right derived functor of $H_Z^0(-)$.
\end{dfn}

If $\mathfrak b$ is an ideal,
  then $Z = V(\mathfrak b)$ is stable under specialization
  and $H_Z^p(-)$ coincides with the ordinary local cohomology functor
  $H_{\mathfrak b}^p(-)$.

Let $Z$ be a subset of $\spec A$
  which is stable under specialization.
If $\mathfrak b$, $\mathfrak b'$ are ideals
  such that $V(\mathfrak b)$, $V(\mathfrak b') \subset Z$,
  then $V(\mathfrak b \cap \mathfrak b') \subset Z$.
Therefore the set $\mathcal F$ of all ideals~$\mathfrak b$ such that
  $V(\mathfrak b) \subset Z$
  is a directed set with respect to the opposite inclusion.
If $\mathfrak b$, $\mathfrak b' \in \mathcal F$
  such that $\mathfrak b' \subset \mathfrak b$,
  then there is a natural transformation
  $\ext_A^p(A/\mathfrak b, -) \to \ext_A^p(A/\mathfrak b', -)$.
Since $H_Z^0(-) =
  \injlim_{\mathfrak b \in \mathcal F} \hom(A/\mathfrak b, -)$,
  we obtain the natural isomorphism
    \begin{equation} \label{eqn:2.1.1}
      H_Z^p(-) =
      \injlim_{\mathfrak b \in \mathcal F} \ext_A^p(A/\mathfrak b, -).
    \end{equation}

The following lemma was essentially given
  by Raghavan~\cite[p.~491]{MR94e:13033}.

\begin{lem}
  Let $M$ be a finitely generated $A$-module.
  Then $\mathcal L
      = \{H_Z^0(M) \mid Z \subset \spec A$ is stable under specialization$\}$
    is a finite set.
\end{lem}

\begin{proof}
  Let $\ass M = \{\mathfrak p_1, \dots, \mathfrak p_r\}$ and
    $0 = M_1 \cap \dots \cap M_r$ be a primary decomposition
    of $0$ in~$M$
    where $\ass M/M_i = \{\mathfrak p_i\}$ for all~$i$.
  Then $H_Z^0(M) = \bigcup_{V(\mathfrak b) \subset Z} 0 \qtn_M \mathfrak b
    = \bigcap_{\mathfrak p_i \notin Z} M_i$.
  Therefore $\# \mathcal L \leq 2^r$.
\end{proof}

We need Cousin complexes to prove Theorem~\ref{thm:1.1}.

Let $M$ be a finitely generated $A$-module.
For a prime ideal $\mathfrak p \in \supp M$,
  the $M$-height of $\mathfrak p$ is defined to be
  $\height_M \mathfrak p = \dim M_{\mathfrak p}$.
If $\mathfrak b$ is an ideal in~$A$ such that $M \ne \mathfrak bM$,
  then let $\height_M \mathfrak b
  = \inf \{\height_M \mathfrak p \mid
  \mathfrak p \in \supp M \cap V(\mathfrak b)\}$.

\begin{dfn}[\cite{MR41:8400}]
  The Cousin complex $(M^\bullet, d_M^\bullet)$ of $M$
    is defined as follows:

  Let $M^{-2} = 0$, $M^{-1} = M$ and
    $d_M^{-2} \colon M^{-2} \to M^{-1}$ be the zero map.
  If $p \geq 0$ and $d_M^{p-2} \colon M^{p-2} \to M^{p-1}$ is given,
    then we put
      $$
        M^p = \bigoplus_
          {
            \substack{
              \mathfrak p \in \supp M \\
              \height_M \mathfrak p = p
            }
          }
        (\coker d_M^{p-2})_{\mathfrak p}.
      $$
  If $\xi \in M^{p-1}$ and
    $\bar\xi$ is the image of $\xi$ in $\coker d_M^{p-2}$,
    then the component of $d_M^p(\xi)$
    in $(\coker d_M^{p-2})_{\mathfrak p}$
    is $\bar\xi/1$.
\end{dfn}

The following theorem contains \cite[Theorems 11.4 and 11.5]{MR91g:13010}.
\begin{thm} \label{thm:2.4}
  Assume that $A$ satisfies \textup{(C1)}--\textup{(C3)} and
    let $M$ be a finitely generated $A$-module satisfying
      \begin{itemize}
        \item[(QU)] $\height \mathfrak p/\mathfrak q + \height_M \mathfrak q
          = \height_M \mathfrak p$
          for any $\mathfrak p$, $\mathfrak q \in \supp M$
          such that $\mathfrak p \supset \mathfrak q$.
      \end{itemize}
  Then there is an ideal $\mathfrak a$ in~$A$
    satisfying the following properties\textup:

  \textup{(1)}~$V(\mathfrak a)$ is the non-Cohen-Macaulay locus of~$M$.
  In particular, $\height_M \mathfrak a > 0$.

  \textup{(2)}~Let $Z$ be a subset of $\spec A$
    which are stable under specialization
    and $n$ an integer.
  If $\height_M \mathfrak p \geq n$ for any $\mathfrak p \in Z \cap \supp M$,
    then $\mathfrak aH_Z^p(M) = 0$ for each $p < n$.

  \textup{(3)} Let $x_1$, \dots, $x_n \in A$ be a sequence.
  If $\height_M (x_1, \dots, x_n)A \geq n$,
    then $\mathfrak a$ annihilates the Koszul cohomology module
    $H^p(x_1, \dots, x_n; M)$ of~$M$
    with respect to $x_1$, \dots, $x_n$
    for any $p < n$.
\end{thm}

\begin{proof}
  Let $M^\bullet$ be the Cousin complex of $M$
    and $\mathfrak a$ the product of all the annihilators
    of all the non-zero cohomologies of~$M^\bullet$.
  Then it is well-defined and satisfies~(1).
  See \cite[Corollary 6.4]{Kawasaki::finiteness}.

  We prove (2).
  Because of \eqref{eqn:2.1.1},
    it is enough to show that $\mathfrak a\ext^p(A/\mathfrak b, M) = 0$
    for any ideal $\mathfrak b$ such that $V(\mathfrak b) \subset Z$
    and for any $p < n$.
  Let $\mathfrak b$ be such an ideal
    and $F_\bullet$ a free resolution of $A/\mathfrak b$.
  The double complex $\hom(F_\bullet, M^\bullet)$ gives two spectral sequences
      \begin{align*}
        {}' E_2^{pq} & = \ext^p(A/\mathfrak b, H^q(M^\bullet)) \Rightarrow
          H^{p+q}(\hom(F_\bullet, M^\bullet)),
      \\
        {}'' E_2^{pq} & = H^p(\ext^q(A/\mathfrak b, M^\bullet)) \Rightarrow
          H^{p+q}(\hom(F_\bullet, M^\bullet)).
      \end{align*}
  The first spectral sequence tells us that
    $\mathfrak a H^k(\hom(F_\bullet, M^\bullet)) = 0$
    for any~$k$.

  On the other hand, ${}'' E_2^{pq} = 0$ if $p < -1$ or if $q < 0$.
  Let $0 \leq p < n$ be an integer
    and $\mathfrak p \in \supp M$ such that $\height_M \mathfrak p = p$.
  Since $\mathfrak b \not\subset \mathfrak p$,
    we find that $\hom(F_\bullet, (\coker d_M^{p-2})_{\mathfrak p})$
    is exact.
  Hence $\hom(F_\bullet, M^p)$ is also exact.
  Thus ${}'' E_2^{pq} = 0$ if $0 \leq p < n$ and
    ${}'' E_2^{-1,q} = \ext^q(A/\mathfrak b, M)$.
  If $k < n$,
    then ${}'' E_2^{p,k-p-1} = {}'' E_2^{p,k-p} = 0$
    whenever $p \ne -1$.
  Therefore $H^{k-1}(\hom(F_\bullet, M^\bullet))
    = {}'' E_2^{-1,k} = \ext^k(A/\mathfrak b, M)$
    is annihilated by $\mathfrak a$.

  Next we consider (3).
  Let $K_\bullet$ be the Koszul complex of $A$
    with respect to $x_1$, \dots,~$x_n$.
  By considering the double complex $\hom(K_\bullet, M^\bullet)$,
    instead of $\hom(F_\bullet, M^\bullet)$,
    we obtain the assertion.
\end{proof}

\section{The proof of Theorem~\ref{thm:1.1}}
Before the proof of Theorem~\ref{thm:1.1},
  we fix some notation.
Let $\mathfrak X$ be the free Abelian group with basis $\spec A$
  and $\mathfrak X_+ = \{\sum k_{\mathfrak p} \mathfrak p \mid
  k_{\mathfrak p} \geq 0$ for all $\mathfrak p\}$.
If $\alpha = k_1 \mathfrak p_1 + \dots + k_n \mathfrak p_n$
  and $\beta = l_1 \mathfrak p_1 + \dots + l_n \mathfrak p_n$
  where $\mathfrak p_i \ne \mathfrak p_j$
  whenever $i \ne j$,
  then we put
    $$
      \alpha \vee \beta = \sum_{i=1}^n \max\{k_i, l_i\} \mathfrak p_i.
    $$
It is clear that $(\alpha \vee \beta) + \gamma
  = (\alpha + \gamma) \vee (\beta + \gamma)$.
Let $\alpha = k_1 \mathfrak p_1 + \dots + k_n \mathfrak p_n \in \mathfrak X_+$
  and $Y$ be a subset of $\spec A$
  which is stable under specialization.
Then we put $\mathfrak b(\alpha, Y)
  = \prod_{\mathfrak p_i \in Y} \mathfrak p_i^{k_i}$.
Since $V(\mathfrak b(\alpha, Y)) \subset Y$,
  Theorem~\ref{thm:1.1} is contained in the following
\begin{thm} \label{thm:3.1}
  Assume that $A$ satisfies \textup{(C1)}--\textup{(C3)}.
  If $M$ is a finitely generated $A$-module,
    then there is $\alpha(M) \in \mathfrak X_+$
    satisfying the following property\textup:

  Let $Y$, $Z$ be subsets of $\spec A$ which are stable under specialization
    and $n$ an integer.
  If
      \begin{itemize}
        \item[(A)]
          $\height \mathfrak p/\mathfrak q + \depth M_{\mathfrak q} \geq n$
          for any $\mathfrak q \in \spec A \setminus Y$ and
          $\mathfrak p \in V(\mathfrak q) \cap Z$,
      \end{itemize}
    then
      \begin{itemize}
        \item[(B)] $\mathfrak b(\alpha(M), Y)$ annihilates
          $H_Z^0(M)$, \dots, $H_Z^{n-1}(M)$.
      \end{itemize}
\end{thm}

We prove this theorem by the Noetherian induction on $\supp M$
  and the induction on the number of associated primes of~$M$.

If $M = 0$, then $\alpha(M) = 0$ obviously satisfies the assertion.
Assume that $M \ne 0$ and that,
  for any finitely generated $A$-module~$M'$,
  there is $\alpha(M')$ satisfying the assertion of Theorem~\ref{thm:3.1}
  if $\supp M' \subsetneq \supp M$ or
  if $\supp M' = \supp M$ and $\# \ass M' < \# \ass M$.
We first prove the following claim.

\begin{claim}
  There is $\alpha'(M) \in \mathfrak X_+$
    satisfying the following property\textup:

  Let $Y$, $Z$ be subsets of $\spec A$ which are stable under specialization
    and $n$ an integer.
  If $Y \cap \ass M = \emptyset$
    and \textup{(A)} holds,
    then \textup{(B)} also does.
\end{claim}

\begin{proof}
  Let $\ass M = \{P_1, \dots, P_r\}$.
  We may assume that $P_1 \not\subset P_2$, \dots, $P_r$
    without loss of generality.
  There is an exact sequence
      $$
        0 \to L \to M \to N \to 0
      $$
    such that $\ass L = \{P_2, \dots, P_r\}$ and
    $\ass N = \{P_1\}$.
  Since $A$ is universally catenary and
    $N$ has the unique minimal prime,
    $N$ satisfies (QU).
  Let $\mathfrak a$ be the ideal obtained by applying Theorem~\ref{thm:2.4}
    to~$N$.
  Then $P_1 \subsetneq \mathfrak a$.
  Since $P_1 \not\subset P_2$, \dots, $P_r$,
    we find that $\mathfrak a \not\subset P_2$, \dots, $P_r$.
  Let $x'' \in \mathfrak a \setminus (P_1 \cup \dots \cup P_r)$.

  Since $\supp L \subsetneq \supp M$ or
    since $\supp L = \supp M$ and $\# \ass L < \# \ass M$,
    there is $\alpha(L) \in \mathfrak X_+$
    satisfying the assertion of Theorem~\ref{thm:3.1}.
  Let $\alpha(L) = k_1 Q_1 + \dots + k_s Q_s$.
  We may assume that $Q_1$, \dots, $Q_{s_0} \not\subset P_1 \cup \dots \cup P_r$
    and $Q_{s_0+1}$, \dots, $Q_s \subset P_1 \cup \dots \cup P_r$.
  Let $x' \in Q_1^{k_1} \cdots Q_{s_0}^{k_{s_0}}
    \setminus P_1 \cup \dots \cup P_r$
    and $x = x' x''$.

  Since $x$ is an $M$-non zero divisor,
    $\supp M/xM \subsetneq \supp M$.
  We want to show that $\alpha'(M) = \alpha(M/xM)$ satisfies
    the assertion of the claim.

  Let $Y$, $Z$ be subsets of $\spec A$
    which are stable under specialization
    and $n$ an integer.
  Assume that $Y \cap \ass M = \emptyset$ and
    $\height \mathfrak p/\mathfrak q + \depth M_{\mathfrak q} \geq n$
    for any $\mathfrak q \in \spec A \setminus Y$
    and $\mathfrak p \in V(\mathfrak q) \cap Z$.
  If $\mathfrak p \in Z \cap \supp N$,
    then $\height \mathfrak p/P_1 + \depth M_{P_1} \geq n$
    because $\supp N = V(P_1)$ and
    $P_1 \notin Y$.
  Since $\depth M_{P_1} = 0$,
    we have
      \begin{equation} \label{eqn:3.1.1}
        \height_N \mathfrak p = \height \mathfrak p/P_1 \geq n
        \quad
        \text{for any $\mathfrak p \in Z \cap \supp N$}.
      \end{equation}
  By using Theorem~\ref{thm:2.4} (2),
    we find that $x'' H_Z^p(N) = 0$ for any $p < n$.

  Let $\mathfrak q \in \spec A \setminus (Y \cup V(x'' A))$
    and $\mathfrak p \in V(\mathfrak q) \cap Z$.
  Since $x'' \notin \mathfrak q$,
    $N_{\mathfrak q}$ is Cohen-Macaulay.
  If $N_{\mathfrak q} \ne 0$,
    then $\mathfrak p \in Z \cap \supp N$ and hence
      \begin{align*}
        \height \mathfrak p/\mathfrak q + \depth N_{\mathfrak q}
        & = \height \mathfrak p/\mathfrak q + \dim N_{\mathfrak q}
      \\
        & = \height_N \mathfrak p \geq n.
      \end{align*}
  Here we used \eqref{eqn:3.1.1}.
  If $N_{\mathfrak q} = 0$,
    then $\depth N_{\mathfrak q} = \infty$ and hence
    $\height \mathfrak p/\mathfrak q + \depth N_{\mathfrak q} \geq n$.
  Since $\mathfrak q \notin Y$,
    the assumption tells us that
    $\height \mathfrak p/\mathfrak q + \depth M_{\mathfrak q} \geq n$.
  Therefore $\height \mathfrak p/\mathfrak q + \depth L_{\mathfrak q} \geq n$.
  Because of the induction hypothesis,
      $$
        \mathfrak b(\alpha(L), Y \cup V(x'' A)) H_Z^p(L) = 0
      $$
    for $p < n$.

  Since $x'' \notin P_1 \cup \dots \cup P_r$,
    $P_1$, \dots, $P_r \notin Y$ and
    $Q_{s_0+1}$, \dots, $Q_s \subset P_1 \cup \dots \cup P_r$,
    we have $Q_{s_0+1}$, \dots, $Q_s \notin Y \cup V(x'' A)$.
  Therefore $x' \in Q_1^{k_1} \cdots Q_{s_0}^{k_{s_0}}
    \subset \mathfrak b(\alpha(L), Y \cup V(x''A))$
    and hence
    $x' H_Z^p(L) = 0$ if $p < n$.
  Since $H_Z^p(L) \to H_Z^p(M) \to H_Z^p(N)$ is exact,
    $x H_Z^p(M) = 0$ if $p < n$.

  Since $x$ is an $M$-non zero divisor,
    $H_Z^0(M) = 0$,
      $$
        0 \to H_Z^{p-1}(M) \to H_Z^{p-1}(M/xM) \to H_Z^p(M) \to 0
      $$
    is exact for $p < n$
    and $\height \mathfrak p/\mathfrak q
    + \depth (M/xM)_{\mathfrak q} \geq n-1$
    for any $\mathfrak q \in \spec A \setminus Y$
    and $\mathfrak p \in V(\mathfrak q) \cap Z$.
  Therefore $\mathfrak b(\alpha'(M), Y) = \mathfrak b(\alpha(M/xM), Y)$
    annihilates $H_Z^p(M)$ if $p < n$.
\end{proof}

Next we give $\alpha(M)$.
Let $\ass M = \{P_1, \dots, P_r\}$ and
  $0 = M_1 \cap \dots \cap M_r$ be a primary decomposition of~$0$ in~$M$
  such that $\ass M/M_i = \{P_i\}$.
Then there are integers $k_1$, \dots, $k_r$
  such that $P_i^{k_i} M \subset M_i$
  for each~$i$.

Let $\{H_Z^0(M) \mid Y \subset \spec A$ is stable under specialization$\}
  = \{L_1, \dots, L_s\}$.
Assume that $L_1 = 0$ and $L_2$, \dots, $L_s \ne 0$.
Since $\supp M/L_i \subsetneq \supp M$ or
  $\supp M/L_i = \supp M$, $\#\ass M/L_i < \# \ass M$,
  there is $\alpha(M/L_i) \in \mathfrak X_+$
  satisfying the assertion of Theorem~\ref{thm:3.1}
  for each $i = 2$, \dots, $s$.
We put $\alpha(M) = \alpha'(M) \vee [\sum k_i P_i + \alpha(M/L_2) \vee \dots
  \vee \alpha(M/L_s)]$.
Then $\alpha(M)$ has required property.

Indeed, let $Y$, $Z$ be subsets of $\spec A$
  which are stable under specialization
  and $n$ an integer.
If $H_Y^0(M) = 0$,
  then $Y \cap \ass M = \emptyset$
  and hence $\mathfrak b(\alpha'(M), Y)$ annihilates
  $H_Z^0(M)$, \dots, $H_Z^{n-1}(M)$.
Assume that $H_Y^0(M) = L_j$ for some $2 \leq j \leq s$.
If $\mathfrak q \in \spec A \setminus Y$ and
  $\mathfrak p \in V(\mathfrak q) \cap Z$,
  then $(L_j)_{\mathfrak q} = 0$ and hence
  $\height \mathfrak p/\mathfrak q + \depth (M/L_j)_{\mathfrak q}
  = \height \mathfrak p/\mathfrak q + \depth M_{\mathfrak q} \geq n$.
Therefore $\mathfrak b(\alpha(M/L_j), Y)$ annihilates
  $H_Z^0(M/L_j)$, \dots, $H_Z^{n-1}(M/L_j)$.
On the other hand, since there is a monomorphism
    $$
      L_j = \bigcap_{P_i \notin Y} M_i
      \hookrightarrow \bigoplus_{P_i \in Y} M/M_i,
    $$
  we find that $\mathfrak b(\sum k_i P_i, Y) L_j = 0$.
Since $H_Z^p(L_j) \to H_Z^p(M) \to H_Z^p(M/L_j)$ is exact,
  $\mathfrak b(\sum k_i P_i + \alpha(M/L_j), Y)$ annihilates
  $H_Z^0(M)$, \dots, $H_Z^{n-1}(M)$.
Thus (B) holds.

If $L_1$, \dots, $L_s$ are all non-zero,
  then we put $\alpha(M) =
  \sum k_i P_i + \alpha(M/L_1) \vee \dots \vee \alpha(M/L_s)$.
We can show that $\alpha(M)$ satisfies the assertion of Theorem~\ref{thm:3.1}
  in the same way as above.
The proof of Theorem~\ref{thm:1.1} is completed.

The following corollary is an improvement of~\cite[Theorem 3.1]{MR94e:13033}.
\begin{cor}
  Assume that $A$ satisfies \textup{(C1)}--\textup{(C3)}.
  If $M$ is a finitely generated $A$-module,
    then there is a positive integer $k$
    satisfying the following property\textup:

  Let $\mathfrak a$, $\mathfrak b$ be ideals in~$A$
    and $n$ an integer.
  If $\height \mathfrak p/\mathfrak q + \depth M_{\mathfrak q} \geq n$
    for any $\mathfrak q \in \spec A \setminus V(\mathfrak b)$
    and $\mathfrak p \in V(\mathfrak a + \mathfrak q)$,
    then $\mathfrak b^k H_{\mathfrak a}^p(M) = 0$ for all $p < n$.
\end{cor}

\begin{proof}
  Let $\alpha(M) = k_1 \mathfrak p_1 + \dots + k_r \mathfrak p_r$
    and $k = k_1 + \dots + k_r$.
  Then $\mathfrak b(\alpha(M), V(\mathfrak b)) \supset \mathfrak b^k$.
\end{proof}

\section{A conjecture of Huneke}

The following theorem is an affirmative answer to Conjecture~2.13
  of~\cite{MR93b:13027}.
Its proof is similar to that of Theorem~\ref{thm:2.4}.

\begin{thm}
  Assume that $A$ satisfies \textup{(C1)}--\textup{(C3)}
    and let $M$ be a finitely generated $A$-module
    satisfying \textup{(QU)}.
  Then there is an ideal $\mathfrak a$ in~$A$
    satisfying the following property\textup:

  \begin{enumerate}
    \item $\height_M \mathfrak a > 0$. \label{item:4.1.1}
    \item Let
          $$
            \begin{CD}
              0 @>>> F^{-n} @>f^{-n}>> F^{-n+1} @>>> \cdots @>>>
              F^{-1} @>f^{-1}>> F^0
            \end{CD}
          $$
        be a complex of finitely generated free $A$-modules
        such that
          \begin{enumerate}
            \item $\rank f^{-n} = \rank F^{-n}$\textup;
            \item $\rank F^i = \rank f^i + \rank f^{i-1}$
              for each $-n < i < 0$\textup;
            \item $\height_M I_{r_i} (f^i) \geq -i$ for each $-n \leq i < 0$
              where $r_i = \rank f_i$ for each~$i$.
          \end{enumerate}
      Then $\mathfrak a H^p(F^\bullet \otimes M) = 0$
        for all $p < 0$.
      Here $I_{r_i}(f^i)$ denotes the ideal generated by
        all the $r_i$-minors of the representation matrix of~$f^i$.
    \end{enumerate}
\end{thm}

\begin{proof}
  Let $M^\bullet$ be the Cousin complex of~$M$
    and $\mathfrak a$ the product of all the annihilators of
    all the non-zero cohomologies of~$M^\bullet$.
  Then $\mathfrak a$ satisfies~(\ref{item:4.1.1}).
  The double complex $F^\bullet \otimes M^\bullet$
    gives a spectral sequence
      $$
        {}' E_2^{pq} = H^p(F^\bullet \otimes H^q(M^\bullet))
        \Rightarrow H^{p+q}(F^\bullet \otimes M^\bullet).
      $$
  It tells us that $\mathfrak a H^p(F^\bullet \otimes M^\bullet) = 0$
    for all~$p$.
  On the other hand, $F^\bullet \otimes M^\bullet$
    gives another spectral sequence
      $
        {}'' E_2^{pq} \Rightarrow H^{p+q}(F^\bullet \otimes M^\bullet)
      $
    where ${}'' E_2^{pq}$ is the cohomology of
      $$
        H^q(F^\bullet \otimes M^{p-1}) \to
        H^q(F^\bullet \otimes M^p) \to
        H^q(F^\bullet \otimes M^{p+1}).
      $$
  If $0 \leq p < n$ and
    $\mathfrak p \in \supp M$ such that $p = \height_M \mathfrak p$,
    then
      $$
        \begin{CD}
          0 @>>> (F^{-n})_{\mathfrak p} @>>> \cdots @>>>
          (F^{-p})_{\mathfrak p}
        \end{CD}
      $$
    is split exact
    and hence $H^q(F^\bullet \otimes M^p) = 0$ if $q < -p$.
  Therefore ${}'' E_2^{pq} = 0$ if $p > 0$ and $p+q < 0$.
  Furthermore ${}'' E_2^{-1,q} = H^q(F^\bullet \otimes M)$ for each $q < 0$.
  Of course, ${}'' E_2^{pq} = 0$ if $p < -1$.
  Thus $H^p(F^\bullet \otimes M) = {}''E_2^{-1,p}
    = H^{p-1}(F^\bullet \otimes M^\bullet)$
    is annihilated by $\mathfrak a$
    if $p < 0$.
\end{proof}

\def\cdprime{$''$} \def\cfac#1{\ifmmode\setbox7\hbox{$\accent"5E#1$}\else
  \setbox7\hbox{\accent"5E#1}\penalty 10000\relax\fi\raise 1\ht7
  \hbox{\lower1.15ex\hbox to 1\wd7{\hss\accent"13\hss}}\penalty 10000
  \hskip-1\wd7\penalty 10000\box7}
  \def\cfgrv#1{\ifmmode\setbox7\hbox{$\accent"5E#1$}\else
  \setbox7\hbox{\accent"5E#1}\penalty 10000\relax\fi\raise 1\ht7
  \hbox{\lower1.05ex\hbox to 1\wd7{\hss\accent"12\hss}}\penalty 10000
  \hskip-1\wd7\penalty 10000\box7}
  \def\cftil#1{\ifmmode\setbox7\hbox{$\accent"5E#1$}\else
  \setbox7\hbox{\accent"5E#1}\penalty 10000\relax\fi\raise 1\ht7
  \hbox{\lower1.15ex\hbox to 1\wd7{\hss\accent"7E\hss}}\penalty 10000
  \hskip-1\wd7\penalty 10000\box7}
  \def\cfudot#1{\ifmmode\setbox7\hbox{$\accent"5E#1$}\else
  \setbox7\hbox{\accent"5E#1}\penalty 10000\relax\fi\raise 1\ht7
  \hbox{\raise.1ex\hbox to 1\wd7{\hss.\hss}}\penalty 10000 \hskip-1\wd7\penalty
  10000\box7} \def\cprime{$'$} \def\Dbar{\leavevmode\lower.6ex\hbox to
  0pt{\hskip-.23ex \accent"16\hss}D}
  \def\dudot#1{\ifmmode{\lineskiplimit=0pt\oalign{$#1$\crcr
  \hidewidth\setbox0=\hbox{\lower1ex\hbox{{\rm\char"7F}}}\dp0=0pt
  \box0\hidewidth}}\relax\else{\lineskiplimit=0pt\oalign{#1\crcr
  \hidewidth\setbox0=\hbox{\lower1ex\hbox{{\rm\char"7F}}}\dp0=0pt
  \box0\hidewidth}}\relax\fi}
  \def\ocirc#1{\ifmmode\setbox0=\hbox{$#1$}\dimen0=\ht0 \advance\dimen0
  by1pt\rlap{\hbox to\wd0{\hss\raise\dimen0
  \hbox{\hskip.2em$\scriptscriptstyle\circ$}\hss}}#1\else {\accent"17 #1}\fi}
  \def\polhk#1{\setbox0=\hbox{#1}{\ooalign{\hidewidth
  \lower1.5ex\hbox{`}\hidewidth\crcr\unhbox0}}}
  \def\rasp{\leavevmode\raise.45ex\hbox{$\rhook$}}
  \def\soft#1{\leavevmode\setbox0=\hbox{h}\dimen7=\ht0\advance \dimen7
  by-1ex\relax\if t#1\relax\rlap{\raise.6\dimen7
  \hbox{\kern.3ex\char'47}}#1\relax\else\if T#1\relax
  \rlap{\raise.5\dimen7\hbox{\kern1.3ex\char'47}}#1\relax \else\if
  d#1\relax\rlap{\raise.5\dimen7\hbox{\kern.9ex \char'47}}#1\relax\else\if
  D#1\relax\rlap{\raise.5\dimen7 \hbox{\kern1.4ex\char'47}}#1\relax\else\if
  l#1\relax \rlap{\raise.5\dimen7\hbox{\kern.4ex\char'47}}#1\relax \else\if
  L#1\relax\rlap{\raise.5\dimen7\hbox{\kern.7ex
  \char'47}}#1\relax\else\message{accent \string\soft \space #1 not
  defined!}#1\relax\fi\fi\fi\fi\fi\fi}
\providecommand{\bysame}{\leavevmode\hbox to3em{\hrulefill}\thinspace}
\providecommand{\MR}{\relax\ifhmode\unskip\space\fi MR }
\providecommand{\MRhref}[2]{%
  \href{http://www.ams.org/mathscinet-getitem?mr=#1}{#2}
}
\providecommand{\href}[2]{#2}

\end{document}